\documentclass[11pt]{amsart}

\usepackage{amsmath,amscd,amssymb}
\usepackage{graphicx}

\newtheorem{theorem}{Theorem}[section]

\newtheorem{proposition}[theorem]{Proposition}

\newtheorem{lemma}[theorem]{Lemma}
\theoremstyle{remark}
\newtheorem{remark}[theorem]{Remark}

\newtheorem{remarks}[theorem]{Remarks}

\newcommand\be{\begin{equation}\label}
\newcommand\ee{\end{equation}}

\newcommand{\N}{\mathbb{N}}
\newcommand{\R}{\mathbb{R}}

\newcommand{\Z}{\mathbb{Z}}

\newcommand{\on}{\operatorname}

\newcommand{\End}{ \on{End} }

\renewcommand{\ker}{ \on{ker}}

\newcommand\dirac{/\kern-1.2ex\partial} 
\newcommand\qu{/\kern-.7ex/} 
\newcommand{\Waff}{W^{\on{a}}} 

\newcommand{\ol}{\overline}

\newcommand\sig{\sigma}
\newcommand\eps{\epsilon}

\renewcommand{\max}{{\on{max}}}

\newcommand{\f}{\frac}

\renewcommand{\l}{\langle}
\renewcommand{\r}{\rangle}

\newcommand{\ti}{\tilde}

\newcommand\vol{\on{vol}}

\newcommand\beqn{\begin{equation}}      
\newcommand\eeqn{\end{equation}}

\newcommand{\mf}{\mathfrak}
\newcommand{\beq}{\begin{eqnarray*}}
\newcommand{\eeq}{\end{eqnarray*}}

\newcommand{\cox}{h^\vee}

\newcommand{\id}{\on{id}}
\newcommand{\reg}{{\on{reg}}}

\newcommand{\ran}{\on{ran}}

\begin{document}
\title[]{Tilings defined by affine Weyl groups}

\vskip.2in
\author{E. Meinrenken}
\address{University of Toronto, Department of Mathematics,
40 St George Street, Toronto, Ontario M4S2E4, Canada }
\email{mein@math.toronto.edu}
%
\begin{abstract}
Let $W$ be a Weyl group, presented as a reflection group on a
Euclidean vector space $V$, and $C\subset V$ an open Weyl chamber. In a
recent paper, Waldspurger proved that the images $(\id-w)(C)$ for 
$w\in W$ are all disjoint, with union the closed cone spanned by 
the positive roots. We prove that similarly, the images $(\id-w)(A)$
of the open Weyl alcove $A$, for $w\in \Waff$ in the affine Weyl
group, are disjoint and their union is $V$.
\end{abstract}

\maketitle 
\setcounter{tocdepth}{2}


\section{Introduction}

Let $W$ be the Weyl group of a simple Lie algebra, presented as a
crystallographic reflection group in a finite-dimensional Euclidean
vector space $(V,\l\cdot,\cdot\r)$.  Choose a fundamental Weyl chamber $C\subset V$, and
let $D$ be its dual cone, i.e. the open cone spanned by the
corresponding positive roots. In his recent paper \cite{wal:rem},
Waldspurger proved the following remarkable result. Consider the
linear transformations $(\id-w)\colon V\to V$ defined by elements
$w\in W$.
\begin{theorem}[Waldspurger] \label{th:wal}
The images $D_w:=(\id-w)(C),\ w\in W$ are all disjoint, and their union is the
closed cone spanned by the positive roots: 
\[ \ol{D}=\bigcup_{w\in W} D_w.\]
\end{theorem}
For instance, the identity transformation $w=\id$ corresponds to $D_{\id}=\{0\}$
in this decomposition, while the reflection $s_\alpha$ defined by a
positive root $\alpha$ corresponds to the open half-line $D_{s_\alpha}=\R_{>0}\cdot\alpha$.  

The aim of this note is to prove a similar result for the \emph{affine}
Weyl group $\Waff$. Recall that $\Waff=\Lambda\rtimes W$ where the
co-root lattice 
$\Lambda\subset V$ acts by translations. Let $A\subset C$ be the Weyl alcove, with $0\in \ol{A}$. 
\begin{theorem}\label{th:affine}
The images $V_w=(\id-w)(A),\ w\in \Waff$ are all disjoint, and their
union is $V$: 
\[ V=\bigcup_{w\in \Waff} V_w.\] 
\end{theorem}
Figure 1 is a picture of
the resulting tiling of $V$ for the root system ${\bf{G}}_2$. 
\begin{figure}\label{fig:g2tile}
\includegraphics[width=10cm]{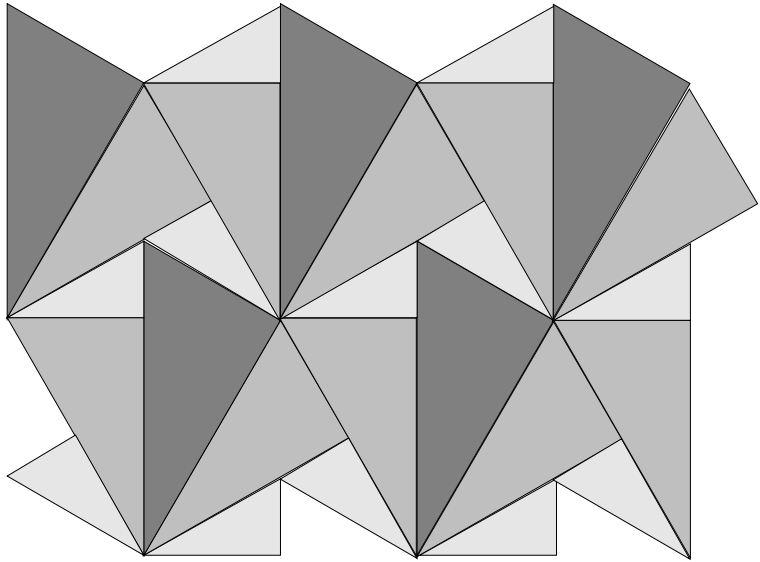}
\caption{The tiling for the root system ${\bf{G}_2}$}
\end{figure}
Up to translation by
elements of the lattice $\Lambda$, there
are five 2-dimensional tiles, corresponding to the five Weyl group
elements with trivial fixed point set.  Letting $s_1,s_2$ denote the
simple reflections, the lightly shaded polytopes are labeled by the
Coxeter elements $s_1s_2,\ s_2s_1$, the medium shaded polytopes by
$(s_1s_2)^2,\ (s_2s_1)^2$, and the darkly shaded polytope by the
longest Weyl group element $w_0=(s_1s_2)^3$.

One also has the following related statement. 
\begin{theorem}\label{th:tau}
  Suppose $S\in \End(V)$ with $||S||<1$. Then the sets
  $V_w^{(S)}=(S-w)(A),\ w\in\Waff$ are all disjoint, and their
  closures cover $V$:
\[ V=\bigcup_{w\in \Waff}\ol{V}_w^{(S)}.\]
\end{theorem}
Note that for $S=0$ the resulting decomposition of
$V$ is just the Stiefel diagram, while for $S=\tau\on{id}$ with
$\tau\to 1$ one recovers the decomposition from Theorem
\ref{th:affine}.

The proof of Theorem \ref{th:affine} is in large parts parallel to
Waldspurger's \cite{wal:rem} proof of Theorem \ref{th:wal}. We will
nevertheless give full details in order to make the paper
self-contained.

\vskip.4in 
\noindent{\bf Acknowledgments:} I would like to thank Bert Kostant
for telling me about Waldspurger's result, and the referee for helpful
comments. I also acknowledge support from an NSERC Discovery Grant and
a Steacie Fellowship.

\section{Notation}
With no loss of generality we will take $W$ to be irreducible.  Let
$\mf{R}\subset V$ be the set of roots,
$\{\alpha_1,\ldots,\alpha_l\}\subset \mf{R}$ a set of simple roots,
and
\[ C=\{x|\ \l\alpha_i,x\r>0,\ i=1,\ldots,l\}\] 
the corresponding Weyl chamber. We denote by $\alpha_\max\in\mf{R}$ the highest
root, and $\alpha_0=-\alpha_\max$ the lowest root.  The open Weyl
alcove is the $l$-dimensional simplex defined as
\[ A=\{x|\ \l\alpha_i,x\r+\delta_{i,0}>0,\ i=0,\ldots,l\}.\]
Its faces are indexed by the proper subsets $I\subset \{0,\ldots,l\}$,
where $A_I$ is given by inequalities $ \l\alpha_i,x\r+\delta_{i,0}>0$
for $i\not\in I$ and equalities $\l\alpha_i,x\r+\delta_{i,0}=0$ for
$i\in I$.  Each $A_I$ has codimension $|I|$. In particular,
$A_i=A_{\{i\}}$ are the codimension $1$ faces, with $\alpha_i$ as
inward-pointing normal vectors. Let $s_i$ be the affine reflections
across the affine hyperplanes supporting $A_i$, 
\[ s_i\colon x\mapsto x-(\l\alpha_i,x\r+\delta_{i,0})\alpha_i^\vee,\ \
\ i=0,\ldots,l,\]
where $\alpha_i^\vee=2\alpha_i/ \l\alpha_i,\alpha_i\r$ is the simple co-root corresponding to $\alpha_i$. The Weyl group
$W$ is generated by the reflections $s_1,\ldots,s_l$, while the affine
Weyl group $\Waff$ is generated by the affine reflections
$s_0,\ldots,s_l$. The affine Weyl group is a semi-direct product
\[ \Waff=\Lambda\rtimes W\]
where the co-root lattice $\Lambda=\Z[\alpha_1^\vee,\ldots,\alpha_l^\vee]\subset
V$ acts on $V$ by translations. For any $w\in \Waff$, we will denote by
$\ti{w}\in W$ its image under the quotient map $\Waff\to W$, i.e.
$\ti{w}(x)=w(x)-w(0)$, and by $\lambda_w=w(0)\in \Lambda$ the
corresponding lattice vector.

The stabilizer of any given element of $A_I$ is the subgroup $\Waff_I$
generated by $s_i,\ i\in I$. It is a finite subgroup of $\Waff$, and
the map $w\mapsto \ti{w}$ induces an isomorphism onto the subgroup
$W_I$ generated by $\ti{s}_i,\ i\in I$. Recall that $W_I$ is itself a Weyl
group (not necessarily irreducible): its Dynkin diagram is obtained from the 
extended Dynkin diagram of the root system $\mf{R}$ by removing all vertices 
that are in $I$.

\section{The top-dimensional polytopes}
For any $w\in\Waff$, the subset 
\[ V_w=(\id-w)(A)\]
is the relative interior of a convex polytope in the affine subspace
$\on{ran}(\id-w)$. Let
\[ \Waff_\reg=\{w\in \Waff|\ (\id-w)\mbox{ is invertible}\}\]
and $W_\reg=W\cap \Waff_\reg$, so that $w\in \Waff_\reg\Leftrightarrow
\ti{w}\in W_\reg$. The top dimensional polytopes $V_w$ are those
indexed by $w\in\Waff_\reg$, and the faces of these polytopes are
$V_{w,I}:=(\id-w)(A_I)$.  For $w\in W_\reg$ and $i=0,\ldots,l$ let
\[ n_{w,i}:=(\id-\ti{w}^{-1})^{-1}(\alpha_i).\]
\begin{lemma}\label{lem:inward}
For all $w\in \Waff_\reg$, the open polytope $V_w$ is given by the inequalities
\[ \l n_{w,i},\xi+\lambda_w\r+\delta_{i,0}>0\]
for $i=0,\ldots,l$. The face $V_{w,I}=(\id-w)(A_I)$ is obtained by
replacing the inequalities for $i\in I$ by equalities. 
\end{lemma}
\begin{proof}
For any $\xi=(\id-w)x\in V$, we have
\[ \l\alpha_i,x\r=\l (\id-\ti{w}^{-1})^{-1}\alpha_i,(\id-\ti{w})x\r=\l
n_{w,i},(\id-\ti{w})x\r
=\l n_{w,i},\xi+\lambda_w\r,\]
since $\ti{w}^{-1}$ is the transpose of $\ti{w}$ under
the inner product $\l\cdot,\cdot\r$. This gives the description of $V_w$
and of its faces $V_{w,I}$. 
\end{proof}
\begin{lemma}\label{lem:pos}
  Suppose $w\in \Waff_\reg,\ i\in \{0,\ldots,l\}$. Then
  \[V_{w,i}=V_{\sigma,i}\subset \ran(\id-\sigma)\] with $\sigma=ws_i$.
  In particular, $\sigma$ is an affine reflection, and $n_{w,i}$ is a
  normal vector to the affine hyperplane $\ran(\id-\sig)$. One has $\l
  n_{w,i},\alpha_i^\vee\r=1$.
\end{lemma}
\begin{proof}
  For any orthogonal transformation $g\in \on{O}(V)$ and any
  reflection $s\in \on{O}(V)$, the dimension of the fixed point set of
  the orthogonal transformations $g,\ gs$ differ by $\pm 1$. Since
  $\ti{w}$ fixes only the origin, it follows that $\ti{\sig}$ has a
  1-dimensional fixed point set. Hence $\on{ran}(\id-\sigma)$ is an
  affine hyperplane, and $\sigma$ is the affine reflection across that
  hyperplane.  Since $s_i$ fixes $A_i$, we have
  $V_{w,i}=(\id-w)(A_i)=(\id-ws_i)(A_i)=V_{\sig,i}\subset
  \ran(\id-\sig)$. By definition
  $n_{w,i}-\ti{w}^{-1}n_{w,i}=\alpha_i$. Hence
\[ -2\l n_{w,i},\alpha_i\r+\l\alpha_i,\alpha_i\r
=||n_{w,i}-\alpha_i||^2-||n_{w,i}||^2=||\ti{w}^{-1}n_{w,i}||^2-||n_{w,i}||^2=0.\qedhere\]
\end{proof}
The following Proposition indicates how the top-dimensional polytopes
$V_{w,i}$ are glued along the polytopes of codimension $1$.
\begin{proposition}\label{prop:first}
Let $\sig\in \Waff$ be an affine reflection, i.e. $\ran(\id-\sig)$ is an
affine hyperplane. Consider 
\begin{equation}\label{eq:no}
\xi\in V_\sigma\backslash \bigcup_{|I|\ge 2} V_{\sigma,I}.
\end{equation}
Then there are two distinct indices $i,i'\in \{0,\ldots,l\}$ such that
$\xi\in V_{\sigma,i}\cap V_{\sigma,i'}$. Furthermore, $w=\sigma s_i$
and $w'=\sigma s_{i'}$ are both in $\Waff_\reg$, so that 
$V_{w,i}=V_{\sigma,i}$ and $V_{w',i'}=V_{\sigma,i'}$, and the
polytopes $V_w,V_{w'}$ are on opposite sides of the affine hyperplane 
$\ran(\id-\sigma)$. 
\end{proposition}
\begin{proof}
  Let $n$ be a generator of the 1-dimensional subspace
  $\ker(\id-\ti{\sig})$. Then $n$ is a normal vector to
  $\ran(\id-\sig)$. The pre-image $(\id-\sig)^{-1}(\xi)\subset V$ is
  an affine line in the direction of $n$. Since $\xi\in V_\sigma$, this
  line intersects $A$, hence it intersects the boundary
  $\partial\ol{A}$ in exactly two points $x,x'$. By 
  \eqref{eq:no}, $x,x'$ are contained in two distinct codimension $1$
  boundary faces $A_i,\ A_{i'}$. Since $n$ is `inward-pointing' at one
  of the boundary faces, and `outward-pointing' at the other, the
  inner products $\l n,\alpha_i\r,\ \l n,\alpha_{i'}\r$ are both
  non-zero, with opposite signs. Let $w=\sigma s_i$ and $w'=\sigma
  s_{i'}$. We will show that $w\in \Waff_\reg$, i.e. $\ti{w}\in
  W_\reg$ (the proof for $w'$ is similar).  Let $z\in V$ with
  $\ti{w}z=z$. Then $\ti{\sig}^{-1}z=\ti{s}_{i} z$, so
\[ (\id-\ti{\sig}^{-1})(z)=(\id-\ti{s}_{i})(z)=\l\alpha_{i},z\r \alpha_{i}^\vee.\]
The left hand side lies in $\on{ran}(\id-\ti{\sig})$, which is
orthogonal to $n$, while the right hand side is proportional to
$\alpha_i$.  Since $\l n,\alpha_{i}\r\not=0$ this is only possible if
both sides are $0$. Thus $z$ is fixed under $\ti{\sig}$, and hence a
multiple of $n$. On the other hand we have $\l\alpha_{i},z\r=0$, hence
using again that $\l n,\alpha_i\r\not=0$ we obtain $z=0$.  This shows
$\ker(\id-\ti{w})=0$.

As we had seen above, $n_{w,i}$ is a normal vector to
$\ran(\id-\sigma)$, hence it is a multiple of $n$. By Lemma
\ref{lem:pos}, it is a positive multiple if and only if $\l
n,\alpha_i\r>0$. But then $\l n,\alpha_{i'}\r<0$, and so $n_{w',i'}$
is a negative multiple of $n$. This shows that $V_w,V_{w'}$ are on
opposite sides of the hyperplane $\ran(\id-\sigma)$.
\end{proof}

Consider the union over $W\subset\Waff$,
\begin{equation}\label{eq:x}
X:=\bigcup_{w\in W} V_w.\end{equation} 
Thus $\bigcup_{w\in\Waff}V_w=\bigcup_{\lambda\in\Lambda}(\lambda+X)$. 
The statement of Theorem \ref{th:affine} means in particular that $X$
is a fundamental domain for the action of $\Lambda$.  Figures 2 and 3
give pictures of $X$ for the root systems ${\bf{B}_2}$ and
${\bf{G}_2}$. The shaded regions are the top-dimensional polytopes
(i.e. the sets $V_w$ for $\id-w$ invertible), the dark lines are the
1-dimensional polytopes (corresponding to reflections), and the origin
corresponds to $w=\id$.
\begin{figure}\label{fig:b2}
\includegraphics[width=6cm]{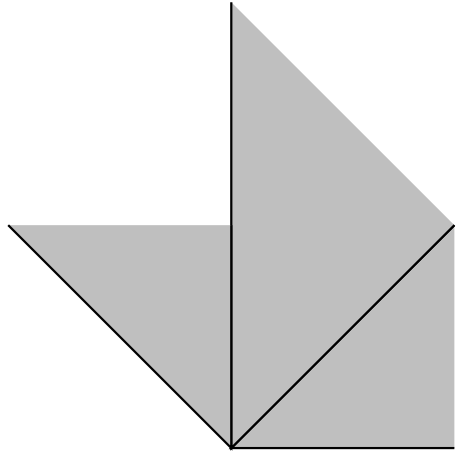}
\caption{The set $X$ for the root system ${\bf{B}_2}$}
\end{figure}
\begin{figure}\label{fig:g2}
\includegraphics[width=6cm]{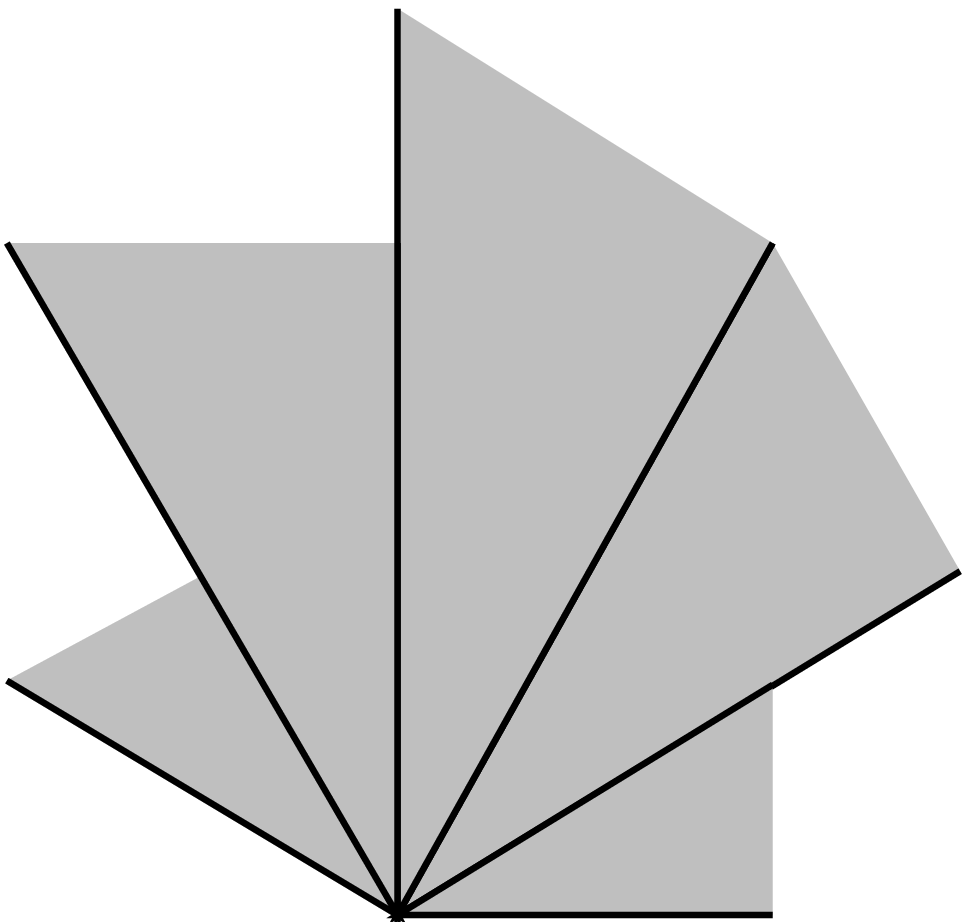}
\caption{The set $X$ for the root system ${\bf{G}_2}$}
\end{figure}

\begin{proposition}\label{prop:cover}
(a) The sets $\lambda+\on{int}(\ol{X}),\ \lambda\in\Lambda$ are
disjoint, and $\bigcup_{\lambda\in\Lambda}\lambda+\ol{X}=V$. 
(b) The open polytopes $V_w$ for $w\in \Waff_\reg$ are disjoint, and 
$\bigcup_{w\in \Waff_\reg}\ol{V}_w=V$. 
\end{proposition}
\begin{proof}
  Since the collection of closed polytopes $\ol{V}_w,\ w\in W_\reg$ is
  locally finite, the union $\bigcup_{w\in \Waff_\reg}\ol{V}_w$ is a
  closed polyhedral subset of $V$.  Proposition \ref{prop:first} shows
  that a point $\xi\in V_{w,i}$ cannot contribute to the boundary of
  this subset unless it lies in
  $\bigcup_{\sigma\in\Waff}\bigcup_{|I|\ge 2}V_{\sigma,I}$. We
  therefore see that the boundary has codimension $\ge 2$, and hence
  is empty since $\bigcup_{w\in \Waff_\reg}\ol{V}_w$ is a closed
  polyhedron.  This proves $\bigcup_{w\in \Waff_\reg}\ol{V}_w=V$, and
  also $\bigcup_{\lambda\in\Lambda}(\lambda+\ol{X})=V$ with $X$ as
  defined in \eqref{eq:x}. Hence the volume $\vol(X)$ (for the
  Riemannian measure on $V$ defined by the inner product) must be at
  least the volume of a fundamental domain for the action of
  $\Lambda$:
\begin{equation}\label{eq:ineq1}
 \vol(X)\ge |W|\vol(A).
\end{equation}
On the other hand, 
$\vol(V_w)=\vol((\id-w)(A))=\det(\id-w)\vol(A)$, so
\begin{equation}\label{eq:ineq2}
\vol(X)\le \sum_{w\in W}\vol(V_w)=\vol(A) \sum_{w\in
  W}\det(\id-w)=|W|\vol(A)
\end{equation} 
where we used the identity \cite[p.134]{bo:li} $\sum_{w\in
  W}\det(\id-w)=|W|$. This confirms $\vol(X)=|W|\vol(A)$.  It follows
that the sets $\lambda+\on{int}(\ol{X})$ are pairwise disjoint, or else the
inquality \eqref{eq:ineq1} would be strict.  Similarly that the sets
$V_w,\ w\in W_\reg$ are disjoint, or else the inequality
\eqref{eq:ineq2} would be strict. (Of course, this also follows from
Waldpurger's Theorem \ref{th:wal} since $C_w\subset D_w$.)  Hence all
$V_w,\ w\in \Waff_\reg$ are disjoint.
\end{proof}

To proceed, we quote the following result from Waldspurger's paper,
where it is stated in greater generality \cite[``Lemme'']{wal:rem}.
\begin{proposition}[Waldspurger] \label{prop:lemme}
  Given $w\in W$ and a proper subset $I\subset \{0,\ldots,l\}$ there exists a
  unique $q\in W_I$ such that
\[ \ker(\id-wq)\cap  \{x\in V|\ \l\alpha_i,x\r>0\mbox{ for all
}i\in I\}\not=\emptyset.\]
\end{proposition}
Following \cite{wal:rem} we use this to prove, 
\begin{proposition}\label{prop:union}
 Every element of $V$ is contained in some $V_w,\
  w\in\Waff$: 
\begin{equation}\label{eq:union} \bigcup_{w\in \Waff}V_w=V.\end{equation}
\end{proposition}
\begin{proof}
  Let $\xi\in V$ be given. Pick $w\in \Waff_\reg$ with $\xi\in
  \ol{V}_w$, and let $I\subset \{0,\ldots,l\}$ with $\xi\in V_{w,I}$.
  Then $x:=(\id-w)^{-1}(\xi)\in A_I$ is fixed under $\Waff_I$. Using
  Proposition \ref{prop:lemme} we may choose $\ti{q}\in W_I$ and
  $n\in V$ such that
\begin{enumerate}
\item $\ti{w}\ti{q}(n)=n$,
\item $\l\alpha_i,n\r>0$ for all $i\in I$
\end{enumerate}
Taking $||n||$ sufficiently small we have $x+n\in A$, and 
\[ (\id-wq)(x+n)=(\id-wq)(x)+(\id-\ti{w}\ti{q})n=(\id-w)(x)=\xi.\]
This shows $\xi\in V_{wq}$.
\end{proof}

\section{Disjointness of the sets $\lambda+X$}
To finish the proof of Theorem \ref{th:affine}, we have to show that
the union \eqref{eq:union} is disjoint. Waldspurger's Theorem
\ref{th:wal} shows that all $D_w=(\id-w)(C),\ w\in W$ are disjoint.
(We refer to his paper for a very simple proof of this fact.) Hence
the same is true for $V_w\subset D_w,\ w\in W$. 
It remains to show that the sets $\lambda+X,\ \lambda\in\Lambda$, with $X$ given by
\eqref{eq:x}, are disjoint. 

The following Lemma shows that the closure $\ol{X}=\bigcup_{w\in
  W}\ol{V}_w$ only involves the top-dimensional polytopes.
\begin{lemma} The closure of the set $X$ is a union over $W_\reg$, 
\[\ol{X}=\bigcup_{w\in W_\reg}\ol{V}_w.\]
Furthermore, $\on{int}(\ol{X})=\on{int}(X)$. 
\end{lemma}
\begin{proof}
  We must show that for any $\xi\in \ol{V}_\sigma,\ \sig \in
  W\backslash W_\reg$, there exists $w\in W_\reg$ such that $\xi\in
  \ol{V}_w$.  Using induction, it is enough to find $\sig'\in W$ such
  that $\xi\in \ol{V}_{\sig'}$ and
  $\dim(\ker(\id-\sig'))=\dim(\ker(\id-\sig))-1$. Let $\pi\colon V\to
  \ker(\id-\sig)^\perp=\ran(\id-\sig)$ denote the orthogonal
  projection. Then $\id-\sig$ restricts to an invertible
  transformation of $\pi(V)$, and $\ol{V}_\sigma$ is the image of
  $\pi(\ol{A})$ under this transformation. We have
\[ \pi(\ol{A})=\pi(\partial\ol{A})=\bigcup_{i=0}^l \pi(\ol{A}_i),\]
and this continues to hold if we remove the index $i=0$ from the right
hand side, as well as all indices $i$ for which $\dim\pi(A_i)<\dim
\pi(V)$. That is, for each point $x\in \pi(\ol{A})$ there exists an
index $i\not=0$ such that $x\in \pi(\ol{A}_i)$, with
$\dim\pi(A_i)=\dim\pi(V)$. Taking $x$ to be the pre-image of $\xi$
under $(\id-\sig)|_{\pi(V)}$, we have $\xi\in \ol{V}_{\sig,i}$ with
$i\not=0$ and $\dim V_{\sig,i}=\dim \ran(\id-\sig)$.  Let $\sig'=\sig
s_i\in W$. Then $V_{\sig,i}=V_{\sig',i}$, hence $
\dim(\ran(\id-\sig'))\ge \dim V_{\sig,i}=\dim(\ran(\id-\sig))$, which
shows $\dim\ker(\id-\sig')\le \dim\ker(\id-\sig)$. By elementary
properties of reflection groups, the dimensions of the fixed point
sets of $\sig,\sig'$ differ by either $+1$ or $-1$. Hence
$\dim(\ker(\id-\sig'))=\dim(\ker(\id-\sig))-1$, proving the first
assertion of the Lemma. 

It follows in particular that the closure of 
$\on{int}(\ol{X})$ equals that of $X$. Suppose $\xi\in \on{int}(\ol{X})$. By Proposition \ref{prop:union}
there exists $\lambda\in \Lambda$ with $\xi\in \lambda+X$. It follows
that $\on{int}(\ol{X})$ meets $\lambda+X$, and hence also
meets $\lambda+\on{int}(\ol{X})$. Since the $\Lambda$-translates of
$\on{int}(\ol{X})$ are pairwise disjoint (see Proposition
\ref{prop:cover}), it follows that $\lambda=0$, i.e.  $\xi \in X$.
This shows $\xi\in X\cap \on{int}(\ol{X})=\on{int}(X)$, hence
$\on{int}(\ol{X})\subset \on{int}(X)$. The opposite inclusion is
obvious.
\end{proof}

Since we already know that the sets $\lambda+\on{int}(X)$ are disjoint, we
are interested in $X\backslash\on{int}(X)\subset \partial
X=\ol{X}\backslash\on{int}(X)$. Let us call a closed codimension $1$
boundary face of the polyhedron $\ol{X}$ `horizontal' if its
supporting hyperplane contains $V_{w,0}$ for some $w\in W_\reg$, and
`vertical' if its supporting hyperplane contains $V_{w,i}$ for some
$w\in W_\reg$ and $i\not=0$. These two cases are exclusive:
\begin{lemma}\label{lem:trans}
  Let $n$ be the inward-pointing normal vector to a codimension $1$
  face of $\ol{X}$.  Then $\l n,\alpha_\max\r\not=0$. In fact, $\l
  n,\alpha_\max\r<0$ for the horizontal faces and $\l
  n,\alpha_\max\r>0$ for the vertical faces.
\end{lemma}
\begin{proof}
  Given a codimension $1$ boundary face of $\ol{X}$, pick any point
  $\xi$ in that boundary face, not lying in $\bigcup_{w\in
    \Waff}\bigcup_{|I|\ge 2} V_{w,I}$. Let $w\in W_\reg$ and $i\in
  \{0,\ldots,l\}$ such that $\xi\in V_{w,i}$, and $n_{w,i}$ is an
  inward-pointing normal vector. By Proposition \ref{prop:first} there
  is a unique $i'\not=i$ such that $\xi\in V_{w',i'}$, where
  $w'=ws_is_{i'}$. Since $V_w,V_{w'}$ lie on opposite sides of the
  affine hyperplane spanned by $V_{w,i}$, and $\xi$ is a boundary
  point of $\ol{X}$, we have $w'\not\in W$. Thus one of $i,i'$ must be
  zero.  If $i=0$ (so that the given boundary face is horizontal) we
  obtain $\l n_{w,0},\alpha_\max\r=-\l n_{w,0},\alpha_0\r<0$.  If
  $i'=0$ we similarly obtain $\l n_{w',0},\alpha_\max\r<0$, hence 
  $\l n_{w,i},\alpha_\max\r>0$.
\end{proof}

\begin{lemma}
  Let $\xi\in X\backslash\on{int}(X)$. Then there exists a vertical
  boundary face of $\ol{X}$ containing $\xi$.  Equivalently, the
  complement $\partial\ol{X} \backslash (X\backslash\on{int}(X))$ is
  contained in the union of horizontal boundary faces.
\end{lemma}
\begin{proof}
  The alcove $A$ is invariant under multiplication by any scalar in $(0,1)$. 
  Hence, the same is true for the sets $V_w$ for $w\in W$, as well as for 
  $X$ and $\on{int}(X)$.  Hence, if
  $\xi\in X\backslash \on{int}(X)$ there exists $t_0>1$
  such that $t\xi\in X\backslash \on{int}(X)$ for $1\le t<t_0$. The
  closed codimension $1$ boundary face containing this line segment is
  necessarily vertical, since a line through the origin intersects the
  affine hyperplane $\{x|\ \l n_{w,0},x-\xi\r=0\}$ in at most one point.
\end{proof}

\begin{proposition}\label{prop:alphamax}
For any $\xi\in X$, there exists $\eps>0$ such that
$\xi+s\alpha_\max\in \on{int}(X)$ for $0< s<\eps$.
\end{proposition}
\begin{proof}
If $\xi\in \on{int}(X)$ there is nothing to show, hence suppose
$\xi\in X\backslash \on{int}(X)$. Suppose first that $\xi$ is not in the union of horizontal
boundary faces of $\ol{X}$. Then there exists an open neighborhood $U$
of $\xi$ such that $U\cap X=U\cap \ol{X}$. All boundary faces of
$\ol{X}$ meeting $\xi$ are vertical, and their inward-pointing normal
vectors $n$ all satisfy $\l n,\alpha_\max\r>0$. Hence,
$\xi+s\alpha_\max\in \on{int}(U\cap \ol{X})=\on{int}(U\cap X)
\subset X$ for $s>0$ sufficiently small. 

For the general case, suppose that for all $\eps>0$, there is $s\in
(0,\eps)$ with $\xi+s\alpha_\max\not\in \on{int}(X)$. We will obtain a
contradiction. Since $\xi$ is contained in some vertical boundary
face, one can choose $t>1$ so that $\xi':=t\xi\in X\backslash \on{int}(X)$, but
$\xi'$ is not in the closure of the union of horizontal boundary faces.
Given $\eps>0$, pick $s\in (0,\eps)$ such that
$\xi+\f{s}{t}\alpha_\max\not\in \on{int}(X)$.  Since $\on{int}(X)$ is invariant under
multiplication by scalars in $(0,1)$, the complement $V\backslash \on{int}(X)$
is invariant under multiplication by scalars in $(1,\infty)$, hence we
obtain $\xi'+s\alpha_\max\not\in \on{int}(X)$.  This contradicts what we have
shown above, and completes the proof.
\end{proof}

\begin{proposition}
The sets $\lambda+X$ for $\lambda\in\Lambda$ are disjoint. 
\end{proposition}
\begin{proof}
Suppose $\xi\in (\lambda+X)\cap (\lambda'+ X)$. By 
Proposition \ref{prop:alphamax}, we can choose $s>0$ so that
  $\xi+s\alpha_\max\in (\lambda+\on{int}(X))\cap 
  (\lambda'+\on{int}(X))$. Since the $\Lambda$-translates of
  $\on{int}(X)$ are disjoint, it follows that $\lambda=\lambda'$.
\end{proof}

This completes the proof of Theorem \ref{th:affine}. We conclude with
some remarks on the properties of the decomposition
$V=\bigcup_{w\in\Waff} V_w$.
\begin{remarks}
\begin{enumerate}
\item The group of symmetries $\tau$ of the extended Dynkin diagram
      (i.e. the outer automorphisms of the corresponding affine Lie
      algebra) acts by symmetries of the decomposition
      $V=\bigcup_{w\in\Waff} V_w$, as follows. Identify the nodes of
      the extended Dynkin diagram with the simple affine reflections
      $s_0,\ldots,s_l$.  Then $\tau$ extends to a group automorphism
      of $\Waff$, taking $s_i$ to $\tau(s_i)$. This automorphism is
      implemented by a unique Euclidean transformation $g\colon V\to
      V$ i.e. $g w g^{-1}=\tau(w)$ for all $w\in \Waff$. Then $g$
      preserves $A$, and consequently
\[ g\,V_w=g(\id-w)(A)=(\id-\tau(w))(A)=V_{\tau(w)},\ \ w\in\Waff.\]
\item It is immediate from the definition that the Euclidean
      transformation $-w\colon V\to V,\ x\mapsto -wx$ takes $V_{w^{-1}}$ into $V_w$:
\[ -w(V_{w^{-1}})=V_w.\]
\item For any positive root $\alpha$, let $s_\alpha$ be the corresponding 
reflection. Then 
$(\id-s_\alpha)(\xi)=\l \alpha,\xi\r\alpha^\vee$, 
where $\alpha^\vee$ is the co-root corresponding to $\alpha$. Hence
$D_{s_\alpha}$ is the relative interior of the line segment from $0$
to $\lambda\alpha^\vee$, where $\lambda$ is the maximum
value of the linear functional $\xi\mapsto \l\alpha,\xi\r$ on the
closed alcove $\ol{A}$. This maximum is achieved at one of
the vertices. Let $\varpi_1^\vee,\ldots,\varpi_l^\vee$ be the
fundamental co-weights, defined by
$\l\alpha_i,\,\varpi_j^\vee\r=\delta_{ij}$ for $i,j=1,\ldots,l$. Let
$c_i\in\N$ be the coefficients of $\alpha_\max$ relative to the simple
roots: $\alpha_\max=\sum_{i=1}^l c_i\alpha_i$. Then the non-zero
vertices of $A$ are $\varpi_i^\vee/c_i$. Similarly let $a_i\in \Z_{\ge
0}$ be the coefficients of $\alpha$, so that $\alpha=\sum_{i=1}^l
a_i\alpha_i$.  Then the value of $\alpha$ at the $i$-th vertex of
$\ol{A}$ is $a_i/c_i$, and $\lambda$ is the maximum of those
values.  Two interesting cases are: (i) If $\alpha=\alpha_\max$, then 
all $a_i/c_i=1$, and $\alpha^\vee=\alpha$. That is, the open line segment 
from the origin to the highest root always appears in the decomposition. 
(ii) If $\alpha=\alpha_i$, then $a_i=1$ while all other $a_j$ vanish. 
In this case, one obtains the open line segment from the origin to 
$\f{1}{c_i}\alpha_i^\vee$.  
\item 
Every $V_w$ contains a distinguished `base point'. Indeed, let
$\rho\in V$ be the half-sum of positive roots, and 
$\cox=1+\l\alpha_\max,\rho\r$ the dual Coxeter number. Then 
$\rho/\cox\in A$, and consequently $\rho/\cox-w(\rho/\cox)\in V_w$. 
\end{enumerate}
\end{remarks}

\section{Proof of Theorem \ref{th:tau}}
The proof is very similar to the proof of Proposition
\ref{prop:cover}, hence we will be brief. 
Each $V_w^{(S)}=(S-w)(A)$ is the interior of a simplex in $V$,
with codimension 1 faces $V_{w,i}^{(S)}=(S-w)(A_i)$. As in the proof of Lemma
\ref{lem:inward}, we see that
\[ n_{w,i}^{(S)}=(S-\ti{w}^{-1})^{-1}\alpha_i\]
is an inward-pointing normal vector to the $i$-th face
$V_{w,i}^{(S)}$. For $S=0$ this simplifies to
\[ n_{w,i}^{(0)}=-w\alpha_i\]
If $w'=w s_i$ we have $V_{w,i}^{(S)}=V_{w',i}^{(S)}$, so that
$n_{w,i}^{(S)}$ and $n_{w',i}^{(S)}$ are proportional. Since 
$n_{w,i}^{(0)}=-n_{w',i}^{(0)}$, it follows by continuity that 
 $n_{w,i}^{(S)}$ is a negative multiple of
  $n_{w',i}^{(S)}$. As a consequence, we see that $V_w^{(S)}, \ V_{w'}^{(S)}$ are on
opposite sides of affine hyperplane supporting $V_{w,i}^{(S)}=V_{w',i}^{(S)}$. Arguing as in the
proof of Proposition \ref{prop:cover}, this shows that 
\[ \bigcup_{w\in\Waff}\ol{V}_w^{(S)}=V.\]
Letting $X^{(S)}=\bigcup_{w\in W} V_w^{(S)}$, it follows that 
$V=\bigcup_{\lambda\in \Lambda}(\lambda+\ol{X}^{(S)})$. Hence  
$\vol({X}^{(S)})\ge |W|\vol(A)$. But
\[ \begin{split}\vol({X}^{(S)})&\le \sum_{w\in W}\vol\big((S-w)(A)\big)\\
  &=\vol(A)\sum_{w\in W}|\det(S-w)|\\
  &=\vol(A)\sum_{w\in W}\det(\id-S w^{-1})=|W|\vol(A),\end{split}\]
using \cite[p.134]{bo:li}. It follows that $\vol({X}^{(S)})=
|W|\vol(A)$, which implies (as in the proof of Proposition
\ref{prop:cover}) that all $\on{int}(\ol{V}_w^{(S)})=V_w^{(S)}$ are
disjoint. This completes the proof.
\begin{remark} Theorem \ref{th:tau}, and its proof, go through for any $S$ in
the component of $0$ in the set 
$\{S\in \End(V)|\ \det(S-w)\not=0\ \forall w\in W\}$. 
For instance, the fact 
that $\det(\id-Sw^{-1})>0$ follows by continuity from $S=0$. 
On the other hand, if e.g. $S$ is a positive matrix with $S>2\id$, the result becomes
false, since then (cf. \cite[p. 134]{bo:li}) $\sum_{w\in W}|\det(S-w)|=\sum_{w\in W}\det(S-w)
=\det(S)|W|$. 
\end{remark}

\def\cprime{$'$} \def\polhk#1{\setbox0=\hbox{#1}{\ooalign{\hidewidth
  \lower1.5ex\hbox{`}\hidewidth\crcr\unhbox0}}} \def\cprime{$'$}
  \def\cprime{$'$} \def\polhk#1{\setbox0=\hbox{#1}{\ooalign{\hidewidth
  \lower1.5ex\hbox{`}\hidewidth\crcr\unhbox0}}} \def\cprime{$'$}
  \def\cprime{$'$}
\providecommand{\bysame}{\leavevmode\hbox to3em{\hrulefill}\thinspace}
\providecommand{\MR}{\relax\ifhmode\unskip\space\fi MR }
\providecommand{\MRhref}[2]{%
  \href{http://www.ams.org/mathscinet-getitem?mr=#1}{#2}
}
\providecommand{\href}[2]{#2}

\end{document}